\documentclass[12pt,a4paper]{article}

	\usepackage{fancyhdr,amsmath, amssymb, amsfonts, dsfont, amsthm, tikz}
    \usepackage{color,epsfig,graphicx}
    \usepackage{color,soul}
    \usepackage{IEEEtrantools}

	\newlength{\guillotine}
	\newtheorem{thm}{Theorem}[section]
	\newtheorem{lem}[thm]{Lemma}
	\newtheorem{notation}[thm]{Notation}
	\newtheorem{defn}[thm]{Definition}
	\newtheorem{example}[thm]{Example}
		\theoremstyle{remark}

	\newcommand{\RR}{\mathbb R}
	\newcommand{\NN}{\mathbb N}

	\newcommand{\calg}{\mathcal G}
	\newcommand{\calc}{\mathcal C}

	\newcommand{\de}{\delta}
	\newcommand{\la}{\lambda}
	\newcommand{\om}{\omega}
	    
	\newcommand{\De}{\Delta}

	\newcommand{\ph}{\phantom}
	
    \newcommand{\un}{\underline}
    \newcommand{\uni}{{\un i}}
    
	\DeclareMathOperator{\Jac}{Jac}
	\DeclareMathOperator{\diam}{diam}
	\DeclareMathOperator{\vol}{vol}

\begin{document}

\title{An elementary proof that the Rauzy gasket is fractal}
\author{Mark Pollicott \& Benedict Sewell\thanks{The  first  author is  partly  supported  by  ERC-Advanced  Grant  833802-Resonances  and EPSRC grant EP/T001674/1, and the second by the Alfr\'ed R\'enyi Young Researcher Fund.}}

\maketitle
\begin{abstract}
	We present an elementary proof that the Rauzy gasket has Hausdorff dimension strictly smaller than two.
\end{abstract}

\begin{figure}[htb]
    \centering
    \includegraphics[width = 0.6\linewidth]{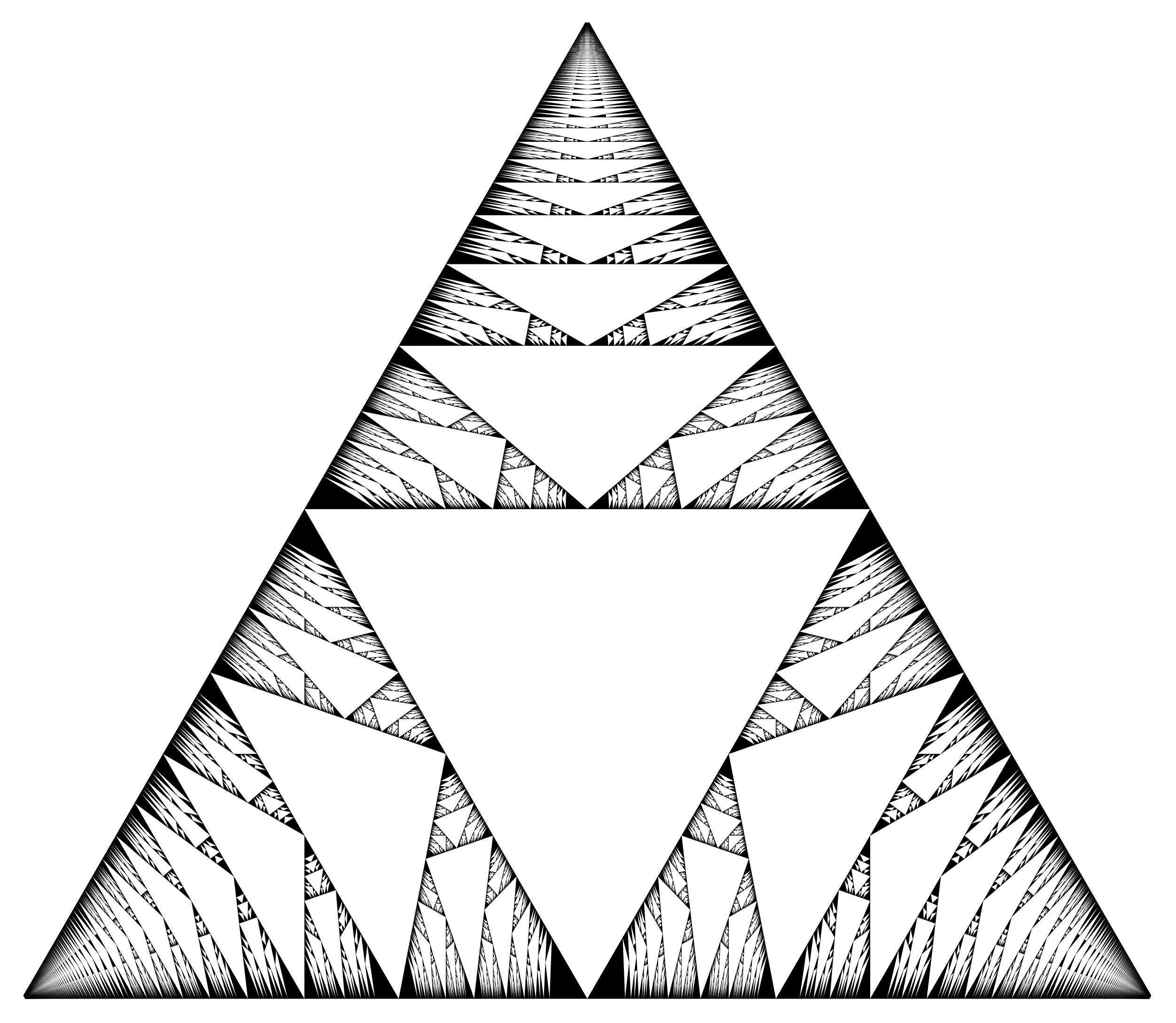}
    \caption{The Rauzy gasket, $\calg_3$.}
\end{figure}

\section{Introduction}

The Rauzy gasket is set which lies within the two-dimensional simplex. It is an important subset of parameter space in numerous dynamical or topological problems (see the introductions of \cite{DD} and \cite{DHS}), having been discovered independently at least three times in different contexts.

It was a conjecture of Arnoux \cite{AS} that the Hausdorff dimension of this set is smaller than two\footnote{It is also a special case of a conjecture of Novikov, that this dimension is strictly between one and two.}, and this was confirmed by \cite{AHS}. In this short note we give a completely elementary proof of this fact.

Our method of proof is fairly flexible, and we illustrate this by applying it to a family of higher-dimensional examples $\calg_d$ ($d\geq 3$) which generalise the Rauzy gasket. That is, the purpose of this note is to give an elementary proof of the following.

\begin{thm}
	$\dim_H(\calg_d) < d-1$ for any $d \geq 3$.
\end{thm}

We note that for the $d=3$ case, this result has been strengthened with explicit (albeit loose) upper bounds \cite{fou,PS}, and complemented with a lower bound \cite{GM}.

In section 2 below, we define $\calg_d$. In section 3 we provide a covering lemma that allows us to reduce to considerations of volumes of certain simplices, and in section 4 we give an explicit formula for these volumes. In section 5, we use the formula to prove some iterative inequalities, which combine with the renewal theorem to get our main technical result. We then verify the assumptions of this technical result in section 6, to complete the proof of Theorem 1.1.

\section{Definition of $\calg_d$}

We begin with the general definition of a Rauzy gasket in the standard $d$-simplex, first defined in \cite{AS}. These were shown to have Hausdorff dimension strictly less than $d$ in \cite{fou}, using estimates of \cite{B}.

\begin{defn}
Fix $d\geq 3$ as the dimension of the ambient space, and let $\De \subset \RR^d$ be the standard $d-1$ simplex:
    $$
        \De = \{x \in[0,1]^d : \|x\| = 1\},
    $$
where $\|x\| = \|x\|_{l^1} = \sum_{j=1}^d v_j$ denotes the usual $l^1$ norm. For $j \in\{1,2,\ldots,d\}$, define the matrix $M_j \in \{0,1\}^{d \times d}$ by 
    $$
            (M_j)_{i,k}
        =
            \begin{cases}
                1 & \text{if $i=j$ or $i=k$,}\\
                0 & \text{otherwise},
            \end{cases}
    $$
i.e., $M_j$ has ones on the main diagonal and the $j$th row, and zeros elsewhere.

We define $\calg_{d-1}$ as the attractor of the projectivised maps $\{T_j\}_{j=1}^d$,
    $$
    		T_j:\De\to \De,
    	\qquad
            T_j (x) 
        = 
            \frac
                {M_j \cdot x}
                {\|M_j\cdot x\|},
    $$
that is, $\calg_{d-1}$ is the unique non-empty compact subset of $\De$ such that
	$$
			\calg_{d-1}
		=
			\bigcup_{j=1}^d
				T_j(\calg_{d-1})
	$$
(see \cite{falc} for equivalent definitions).
\end{defn}

\begin{example}[$d=3$]
	The original Rauzy gasket, depicted in Figure 1, corresponds to $d=3$. The three matrices are
		$$
				M_1 
			=
				\left(
				\begin{array}{ccc}
				 1 & 1 & 1 \\
				 0 & 1 & 0 \\
				 0 & 0 & 1 \\
				\end{array}
				\right),
			\qquad
				M_2
			=
				\left(
				\begin{array}{ccc}
				 1 & 0 & 0 \\
				 1 & 1 & 1 \\
				 0 & 0 & 1 \\
				\end{array}
				\right),
			\qquad
				M_3
			=
				\left(
				\begin{array}{ccc}
				 1 & 0 & 0 \\
				 0 & 1 & 0 \\
				 1 & 1 & 1 \\
				\end{array}
				\right).
		$$
	and the first map is 
		$$
				T_1(x,y,z) 
			= 
				\left(
					\frac 1 {2-x},
					\frac y {2-x},
					\frac z {2-x}
				\right)
%				\frac
%					1
%%					{(1,y,z)}
%					{2-x}
%				{(1,y,z)}
			.
%			\qquad
%				T_2(x,y,z)
%			=
%				\frac{(x,1,z)}{2-y}.
		$$
\end{example}

We now introduce some notation to be used throughout. We shall also write $|\uni| = n$ as shorthand for $\uni \in \{1,2,\ldots,d\}^n$.

\begin{notation}
	Throughout we write, for $\uni = (i_1,\ldots, i_n) \in \{1,2,\ldots, d\}^n$,
	\begin{itemize}	
	    \item $M_\uni := M_{i_1}M_{i_2}\cdots M_{i_n}$,
	    
	    \item $T_\uni := T_{i_1}\circ T_{i_2}\circ \cdots \circ T_{i_n}$,
	 and
	    \item $\De_\uni := T_\uni(\De)$.
	\end{itemize}
\end{notation}

\section{Covering lemma}

Our first step in the proof of Theorem 1.1 is the following lemma. 

\begin{lem}
    Given $\de \in (0,1)$, $\dim_H(\calg_{d}) \leq d + \de - 2$, whenever
        $$
                X_n 
            := 
                \sum_{|\uni| = n}
                    \vol_{d-1}(\De_\uni)^{\de} 
            \to
                0
        $$
    as $n \to \infty$, where $\vol_{d-1}(\De_\uni)$ denotes the $(d-1)$-dimensional volume of $\De_\uni$.
\end{lem}

Note that, where appropriate, we  write $f \lesssim g$ to mean that there exists a constant $C>0$, depending only on $d$ and $\delta$, such that $f \leq Cg$.

\begin{proof}
	From the definition of Hausdorff dimension \cite{falc}, to show $\dim_H(\calg_d) \leq d+\delta - 2$ it suffices to exhibit a family of open covers $\{\calc_n\}_{n=1}^\infty$ of $\calg$ such that
		$$
			\sum_{S \in \calc_n} 
				\diam(S)^{d+\de-2}
			\to 0
		$$
	as $n \to \infty$ (in particular, this implies $\max_{S\in\calc_n} \diam(S) \to 0$). To define these covers, it follows from the definition of $\calg_d$, and the fact that $\calg_d\subset \De$, that
		$$
				\calg_d
			\subset
				\bigcup_{|\uni|=n}
					\De_\uni.
		$$
	Thereby, providing a cover of each $\De_\uni$ by open balls and taking the union gives an open cover of $\calg_n$. For each $|\uni| = n$, the construction is as follows.

	\begin{figure}[tbh]
		\centering
		\includegraphics{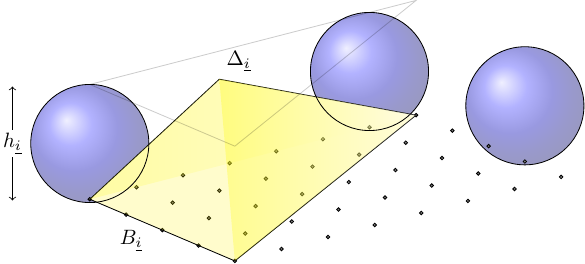}
		\caption{A collection of balls of diameter $h_\uni$ covering a copy of $B_\uni$. The bases of the balls are denoted by dots. Scaling each ball by a factor of $\sqrt 2$ about its centre gives a cover of $\De_\uni$.}
		\label{fig:cover}
	\end{figure}

	Since $T_\uni$ is the projectivisation of an injective linear map, $\De_\uni$ is a $(d-1)$-simplex. Choose a maximal-volume face $B_\uni$ of $\De_\uni$ as its base, so that $h_\uni$, its height measured from $B_\uni$, is the smallest such height. A simple induction shows that $B_\uni$ is contained in a $(d-2)$-dimensional hypercuboid, whose side lengths are at least $h_\uni$ and whose volume is at most $(d-2)!\vol_{d-2}(B_\uni)$. It follows that this cuboid and hence $B_\uni$ can be covered with $n_\uni$ open balls of radius $h_\uni$, where 
		$$
				n_\uni
			\lesssim
%				\left(
					\frac
						{\vol_{d-2}(B_\uni)}
						{h_\uni^{d-2}} 
%				\right)
			\lesssim
%				\left(
					\frac
						{\vol_{d-2}(B_\uni)^{d-1}}
						{\vol_{d-1}(\De_\uni)^{d-2}} 
%				\right)	
			.
		$$

	Treating $\De_\uni \subset B_\uni \times [0,h_\uni]$, taking such a cover of $B_\uni \times \{h_\uni/2\}$ and enlarging the balls by a factor of $\sqrt 2$ gives a cover of $\De_\uni$. Doing this construction simultaneously for all such $|\uni|=n$ defines the cover $\calc_n$. Then, since
		$$
				\sum_{S \in \calc_n}
					\diam(S)^{d+\de-1}
			\lesssim
				\sum_{|i|=n}
					\frac
						{\vol_{d-2}(B_\uni)^{d-1}}
						{\vol_{d-1}(\De_\uni)^{d-2}} 
					\cdot
					\left(
						\frac
							{\vol_{d-1}(\De_\uni)}
							{\vol_{d-2}(B_\uni)}	
					\right)
					^{d+\de-2}
			=
				\sum_{|i|=n}
					\vol_{d-1}(\De_\uni)^{\de}
					\vol_{d-2}(B_\uni)^{1-\de}
		$$
	and $\vol_{d-2}(B_\uni) \lesssim 1,$ our assumption gives the required convergence to zero as $n\to\infty$.
\end{proof}

\section{Volume formula}

We now show that %the Jacobean of $T_\uni$ and 
the volume of $\De_\uni$ can be expressed simply in terms of the entries of $M_\uni$.

\begin{lem}
    For any tuple $\uni$,
        $$
                \frac
                	{\vol_{d-1}(\De_\uni)}
	                {\vol_{d-1}(\De)}
            =
                \nu(M_\uni),
        $$
    where, for any $N \in \RR^{d\times d}$,
        $$
                \nu(N)
            :=
                \prod_{j=1}^d
                \left(
                    \sum_{k=1}^d
                        N_{k,j}
                \right)^{-1}
            = 
                \prod_{j=1}^d
                    \|N \cdot e_j\|^{-1}               
                ,
        $$
    and where $e_j = (0,\ldots,0,1,0,\ldots,0)$ is the $j$th standard basis vector in $\RR^d$, i.e., vertex of $\De$.
\end{lem}

\begin{proof}
	We define  $\De_\uni^\ast := \{\la v : \la \in [0,1],\, v\in \De_\uni\}$, i.e., the $d$-simplex with vertices consisting of those of $\De_\uni$ plus the origin, and define $\De^\ast$ analogously. Integrating appropriately gives
		\begin{equation}
		\label{eq:ratio of volumes}
	            \frac
	          	 	{\vol_{d-1}(\De_\uni)}
		            {\vol_{d-1}(\De)}
			=
	            \frac
	          	 	{\vol_d(\De_\uni^\ast)}
		            {\vol_d(\De^\ast)}.
		\end{equation}
	The two $d$-simplices are related by a linear action: $\De_\uni^\ast = M_\uni^\ast \cdot \De^\ast$, where
		$$
				M_\uni^\ast 
			:= 
				\big(
						T_\uni(e_1) 
					\big| 
						T_\uni(e_2)
					\big| 
					 	\cdots
					\big| 
					 	T_\uni(e_d)
		 	 	\big)
		 	 =
				\left(
						\frac
			                {M_\uni \cdot e_1}
			                {\|M_\uni\cdot e_1\|}
					\;\bigg|\; 
						\frac
			                {M_\uni \cdot e_2}
			                {\|M_\uni\cdot e_2\|}
					\;\bigg| 
					 	\cdots
					\bigg|\;
						\frac
			                {M_\uni \cdot e_d}
			                {\|M_\uni\cdot e_d\|}
		 	 	\right).	 	
		$$
	That is, $M_\uni^\ast$ is $M_\uni$ with each $j$th column multiplied by a factor of $\|M_\uni\cdot e_j\|^{-1}$. Therefore, using that $\det(M_j) = 1$ for each $j=1,\ldots,d$, the right hand side of \eqref{eq:ratio of volumes} is none other than
		$$
			\det(M_\uni^\ast) = \nu(M_\uni) \det(M_\uni) = \nu(M_\uni),
		$$
	as required. 
\end{proof}

%For simplicity, we will assume that $\vol_{d-1}(\De) = 1$ from this point.

\section{Renewal theorem}

We first give a convenient partition of $\{1,\ldots,d\}^n$ and decomposition of $X_n$. Note that these are incomplete, and this gives rise to the remainder term $r_n$ in Lemma 5.3.

\begin{defn}
    For each $n > k \geq 1$, let
        $$
            	A_{n,k} 
            = 
            	\big\{
            			\uni 
            		= 
            			(i_1,i_2,\ldots,i_n)
            		\in 
            			\{1,2,\ldots, d\}^n 
            	\mid
            			 i_1 = \cdots = i_k \neq i_{k+1}
            	\big\}
        $$ 
    and
        $$
                X_{n,k} 
            = 
                \sum_{\uni \in A_{n,k}} 
                    \vol_{d-1}(\De_\uni)^{\de}.
        $$
\end{defn}

The following lemma allows us to focus on the convergence of $X_{n,1}$ in place of $X_n$.

\begin{lem}
	$X_{n} \lesssim X_{n+2,1}$.
%
%    There exists $C = C(\de)>0$ such that $CX_{n} \leq X_{n+2,1}$ for all $n \in \NN$.
\end{lem}

\begin{proof}
	By a change of variables, we simply have 
	$$
			X_{n+2,1} 
		\geq 
			\sum_{|\uni|=n}
				\vol_{d-1}(T_1T_2(\De_\uni))^\de
		\geq
			\min_{x \in \De}
				\big(\Jac_x T_1T_2\big)^\de
			\sum_{|\uni|=n}
				\vol_{d-1}(\De_\uni)^\de	
		=:
			C X_n.
	$$
This constant $C$ is positive as $T_1$ and $T_2$ are injective.
\end{proof}

The next important lemma is our main tool to guarantee that $X_{n,1}$ converges to zero.

\begin{lem} 
	For all $n>k\geq 1$,
	    $$
	            X_{n+1,k+1}
	        \leq 
	            b_k X_{n,k},
	    \qquad
		\text{and}
	    \qquad
	            X_{n+1,1}
	        \leq 
	            \sum_{j=1}^{n-1}
	                a_k 
	                X_{n,j} 
	        +
	            r_n,        
	    $$
	where
		$$
	            b_k 
	        =
	        	\max_{v \in R_k\cap\De_1}
	            \left\|
	            	M_1\cdot v
	            \right\|^{-d\de} ,
		\qquad
	%
%	\text{and}
%	%
%	    \qquad
	        	a_k 
	        = 
	        	\max_{v \in R_k\cap\De_1}
	        	\sum_{j=2}^d
	            \left\|
	            	M_j\cdot v
	            \right\|^{-d\de} 
	    $$
	and $r_n \lesssim n^{(1-d)\de}$, and where
		$$
				R_k
			:=
				\emph{cl}
				\left(
				\bigcup_{j=1}^d
					T_j^k(\De)
					\setminus
					T_j^{k+1}(\De)
				\right),
			\qquad
				i.e.,
			\qquad
				R_k\cap\De_1
			=
				\emph{cl}
					\left(
						T_1^k(\De)
					\setminus
						T_1^{k+1}(\De)
					\right).	
		$$
\end{lem}

\begin{proof}
		First note that, by the definition of $A_{n,k}$, if $\uni = (i_1,\ldots,i_n) \in A_{n,k}$ then
		\begin{itemize}
			\item $(i_1;\uni):= (i_1,i_1,i_2,\ldots,i_n) \in A_{n+1,k+1}$; and
			\item $(j;\uni) := (j,i_1,i_2,\ldots,i_n) \in A_{n+1,1}$ for any $j \neq i_1$.
		\end{itemize}
		A similar statement applies for $\uni \in \{1\}^n \cup \{2\}^n \cup \cdots \cup \{d\}^n$.
	
		From this, $\uni' \in A_{n+1,k+1}$ if and only if there exists a unique $\uni \in A_{n,k}$ such that $\uni' = (i_1;\uni)$. Consequently,
			\begin{equation}
					X_{n+1,k+1} 
				=
					\sum_{\uni' \in A_{n+1,k+1}} 
						\vol_{d-1}(\De_{\uni'})^\de
				=
					\sum_{i \in A_{n,k}}
						\vol_{d-1}(T_{i_1}(\De_\uni))^\de.
				\label{eqr-X n+1 k+1 over X n k}
			\end{equation}
		Using the formula in Lemma 4.1 and that
			$$
					\uni \in A_{n,k} 
				\implies
					\De_\uni
				\in 
					R_k
				:=
					\text{cl}
					\left(
					\bigcup_{j=1}^d
						T_j^k(\De)
						\setminus
						T_j^{k+1}(\De)
					\right)
			$$
		(i.e., $T_\uni(e_j) \in R_k \cap \De_{\uni_1}$ for each $j$), we have
			\begin{IEEEeqnarray}{rClCr}
					\frac
						{\vol_{d-1}(T_{i_1}(\De_\uni))^\de}
						{\vol_{d-1}(\De_{\uni'})^\de}
				&=& 
					\prod_{j=1}^d
						\frac
							{\|M_{i_1}M_\uni \cdot e_j\|^{-\de}}
							{\|M_\uni \cdot e_j\|^{-\de}}
%				\nonumber 
%	 			\\
				&=& 
					\prod_{j=1}^d		
						\left\|
							M_{i_1} \cdot 
							\frac
								{M_\uni \cdot e_j}
								{\|M_\uni \cdot e_j\|}
						\right\|^{-\de}
				\nonumber 
	 			\\
				&=& 
					\prod_{j=1}^d		
						\|M_{i_1} \cdot T_\uni (e_j)\|^{-\de}
%				\nonumber 
%	 			\\
				&\leq&
					\prod_{j=1}^d	
						\max_{v \in R_k \cap \De_{i_1}}\|M_{i_1} \cdot v\|^{-\de}
				\nonumber 
	 			\\
				&=&
					\prod_{j=1}^d			
						\max_{v \in R_k \cap \De_{1_{\ph1}}}\|M_{1} \cdot v\|^{-\de}				
%				\nonumber 
%	 			\\
				&=&
%					\ph{\prod_{j=1}{}}
					\max_{v \in R_k \cap \De_{1_{\ph1}}}
						\|M_{1} \cdot v\|^{-d\de}	
				\nonumber
				\\
				&=:& b_k,
				\label{eqr-b_n intro array 2}
	%					.
			\end{IEEEeqnarray}
		using symmetry and linearity.
		Applying this estimate to \eqref{eqr-X n+1 k+1 over X n k} gives the required inequality:
			\begin{align*}
					X_{n+1,k+1}
				&\leq
					\sum_{i \in A_{n,k}}
						b_k\,
						\vol_{d-1}(\De_\uni)^\de
				= 
					b_k X_{n,k}.
			\end{align*}

		The proof of the second inequality in the lemma is slightly more nuanced. From our first consideration, we have
			\begin{equation*}
					X_{n+1,1}
				=
					\sum_{|\uni|=n}
					\sum_{\substack{1\leq \om\leq d:\\\om \neq i_1}}
						\vol_{d-1}(T_\om(\De_\uni))^\de,
%					\label{eqr-X_n+1,1 = F(i)s}			
			\end{equation*}
		and we bound the internal sum in two cases:
		\\[8pt]
		\noindent\textbf{Case 1:} $\uni \in A_{n,k}$ \textit{for some} $k$. This case is similar to the proof of \eqref{eqr-b_n intro array 2}, but this time we also apply the AM-GM inequality:
		\allowdisplaybreaks
			\begin{IEEEeqnarray}{rClCr}
					\frac{\sum_{\om \neq i_1} 
						\vol_{d-1}(T_\om(\De_\uni))^\de}
						{\vol_{d-1}(\De_\uni)^\de}
				&=&
					\sum_{\om\neq i_1}
						\prod_{j=1}^d
						\frac
							{\|M_{\om}M_\uni \cdot e_j\|^{-\de}}
							{\|M_\uni \cdot e_j\|^{-\de}}
				&=&
					\sum_{\om \neq i_1}	
					\prod_{j=1}^d
						\|M_{\om} \cdot T_\uni (e_j)\|^{-\de}
				\nonumber 
	 			\\ % AM-GM inequality
				&\leq&
					\frac 1 d
					\sum_{\om\neq i_1}
						\sum_{j=1}^d			
							\|M_{\om} \cdot T_\uni (e_j)\|^{-d\de}
				&=&
					\frac 1 d
					\sum_{j=1}^d	
						\sum_{\om\neq i_1}		
							\|M_{\om} \cdot T_\uni (e_j)\|^{-d\de}
				\nonumber 
	 			\\
				&\leq&	
					\max_{v \in R_k \cap \De_{i_1}}
						\sum_{\om \neq i_1}\|M_{\om} \cdot v\|^{-d\de}	
				&=&
					\max_{v \in R_k\cap\De_1}
					\sum_{j=2}^d
						\|M_j \cdot v\|^{-d\de}
				\nonumber 
	 			\\
				&=:&
					a_k.
					\label{eqr-a_n intro array 2}
			\end{IEEEeqnarray}
		Summing over $\uni\in A_{n,k}$ hence gives
			$$
					\sum_{i \in A_{n,k}}
					\sum_{\om \neq i_1} 
						\vol_{d-1}(T_\om(\De_\uni))^\de
				\leq		
					\sum_{i \in A_{n,k}}
						a_k
						\vol_{d-1}(\De_\uni)^\de
				=			
					a_k X_{n,k}.
			$$

		\noindent\textbf{Case 2:} $\uni \in \{1\}^n \cup \{2\}^n \cup \cdots \cup \{d\}^n$. This is an explicit calculation, using that all $\binom d 2$ of the summands for this case are equal:
			$$
					r_n
				:=
					\sum_{j=1}^d
					\sum_{\om \neq j} 
						\vol_{d-1}(T_\om T_j^n (\De))^\de
				\lesssim
					\vol_{d-1}(T_1T_2^n(\De))^\de
				=
					\big(
						2^{1-d}
						(2n+1)^{-1}
						(n+1)^{2-d}
					\big)^\de
				\lesssim
					n^{\de(1-d)},
			$$
		using the formula in Lemma 4.1 and the explicit form of $M_1M_2^n$. This case completes the proof of the second inequality and hence of the lemma.
\end{proof}

The only missing piece, before we prove the main technical result of this paper,  is to furnish $a_k$ and $b_k$ with values. The proof is simple, but for convenience we defer it to the appendix.

\begin{lem}
	For each $k\in \NN$,
		$$
	            b_k 
	        = 
	            \left(
	                \frac
	                    {k+2}
	                    {k+3}
          		\right)^{d\de} 
		\qquad
	\text{and}
	    \qquad
	        	a_k 
	        = 
	        	\left(
	        		\frac
	        			{k+1}
	        			{2k+1}
	        	\right)^{d\de} 
	        + 
	        	2^{-d\de}
	        	(d-2).
	    $$
\end{lem}

We now use the proof of a particular case of the renewal theorem \cite[p330]{fell} to conclude the following result from the previous.

\begin{thm}
    If $\de > (d-1)^{-1}$ and 
        $
                \sum_{k=1}^\infty
                    a_k
                    \prod_{j=1}^{k-1}
                        b_j
            < 
                1,
        $
    i.e.,
  		$$
				3^{d\de}
				\sum_{k=1}^\infty
					\left(
						\frac 
							{k+1}
							{2k+1}
					\right)^{d\de} 
					(k+2)^{-d\de}
			+ 
				2^{-d\de}
				(d-2) 
				(k+2)^{-d\de}
			<
				1,
	$$    

    then $X_{n,1} \to 0$ as $n \to \infty$. Consequently, $\dim_H(\calg_d)\leq d-2+\de$
\end{thm}

\begin{proof}
	Applying the first inequality of Lemma 4.3 to the summands of the second gives
		$$
				X_{n+1,1}
			\leq
				\sum_{k=1}^{n-1}
					a_k
					\prod_{j=1}^{k-1}
						b_j\,
					X_{n+1-k,1}
			+ r_n.
		$$
	Writing $\la_k = a_k\prod_{j=1}^{k-1}b_j$ for succinctness, we then have
		\begin{IEEEeqnarray*}{lCcClCl}
				\sum_{n=1}^N
					X_{n+1,1}
			&\leq&
				\sum_{n=1}^N
				\sum_{k=1}^{n-1}
					\la_k
					X_{n+1-k,1}
			+
				\sum_{n=1}^N
					r_n
			&=&
				\sum_{k=1}^{N-1}
					\la_k
				\sum_{n=k+1}^{N}
					X_{n+1-k,1}
			&+&
				\sum_{n=1}^N
					r_n		
			\\
			&&
			&\leq&
				\sum_{k=1}^{N-1}
					\la_k	
				\sum_{j=1}^{N+1}
					X_{j,1}
			&+&
				\sum_{n=1}^N
					r_n,									
		\end{IEEEeqnarray*}
	i.e.,
		$$
				\sum_{n=1}^N
					X_{n+1,1}
			\leq
				\frac
					{ x_1\sum_{k=1}^{N-1} \la_k + \sum_{n=1}^N r_n}
					{1- \sum_{k=1}^N \la_k}.
		$$
	The right hand side is bounded in $N$ by our assumptions, and the result follows.
\end{proof}

\section{Proof of Theorem 1.1}

The remainder of the proof of Theorem 1.1 is to show that Theorem 5.5 holds for $\de=1$ and any $d\geq 3$, since by continuity, it will then apply for any $\de<1$ sufficiently close to 1. More explicitly, we wish to show that
	$$
			\sum_{k=1}^\infty
				\left(
					\frac
						{3(k+1)}
						{(2k+1)(k+2)}
				\right)^d
			+ 
				(d-2)
				\left(
					\frac3{2(k+2)}
				\right)^{d} 
		=
			\sum_{k=1}^\infty
				\left(
						\frac 1{2k+1}
					+
						\frac1{k+2}
				\right)^d
+ 
			(d-2)
				\left(
					\frac3{2(k+2)}
				\right)^{d}			
	$$
is strictly less than 1. This uses elementary calculus.
Since each term on the right hand side is decreasing in $d$ for $d\geq 3$, it suffices to show just the $d=3$ case, i.e., that
	\begin{equation}
		\frac{27}8
			\sum_{k=1}^\infty
				\left(
					\frac
						{2(k+1)}
						{(2k+1)(k+2)}
				\right)^3
			+ 
				\left(
					k+2
				\right)^{-3}
		<
			1.
	\label{eq:d=3 renewal}			
	\end{equation}
This follows by simply bounding the tail of the sum by an integral. For any $n\in\NN$, the left hand side of \eqref{eq:d=3 renewal} is at most
	\begin{gather*}
			\frac{27}8
				\sum_{k=1}^n
					\left(
						\frac
							{2(k+1)}
							{(2k+1)(k+2)}
					\right)^3
				+ 
					\left(
						k+2
					\right)^{-3}
		\quad+\quad
			\frac{27}8
			\int_{n-1}^\infty
				8(2x+1)^{-3} + (x+2)^{-3}
			\;\mathrm dx
		\\
		{}=
			\frac{27}8
			\sum_{k=1}^n
				\left(
					\frac
						{2(k+1)}
						{(2k+1)(k+2)}
				\right)^3
			+ 
				\left(
					k+2
				\right)^{-3}
		\quad+\quad
			\frac {27}{4(2n-1)^{2}} 
		+ 
			\frac {27}{16(n+1)^{2}},		
		\end{gather*}
and for $n=3$, this expression equals
	$$
			\left(\frac23\right)^3
		+
			\left(\frac12\right)^3
		+
			\left(\frac9{20}\right)^3
		+
			\left(\frac38\right)^3
		+
			\left(\frac{12}{35}\right)^3
		+
			\left(\frac3{10}\right)^{3}
		+
			\frac{27}{4 \times 5^2}
		+
			\frac{27}{16 \times 4^2}
		= 
			\frac{574898507}{592704000}
		<1.
	$$
Hence \eqref{eq:d=3 renewal} holds, completing the proof of the theorem.

\appendix

\section{Proof of Lemma 5.4}
Noting that $\|M_j \cdot v\| = (2-v_j)$ and
		$$
				R_k\cap\De_1
			:= 
				\left\{
					v \in \De
				\;
				\bigg|
				\;
					\frac{k}{k+1}
				\leq
					v_1
				\leq
					\frac{k+1}{k+2}	
			\right\},
		$$
we have
	$$
			b_k
		=
			\max
				\left\{
					(2-v_1)^{-d\de}
				\;
				\bigg|
				\;
					\frac{k}{k+1}
				\leq
					v_1
				\leq
					\frac{k+1}{k+2}					
				\right\}
		=
			\left(
				2 - \frac{k+1}{k+2}	
			\right)^{-d\de} 
		=
			\left(
				\frac{k+2}{k+3}
			\right)^{d\de}.
	$$

Regarding the value of $a_k$, we write $a_k = \max_{R_k\cap\De_1} f$, where
	$$
			f(v)
		:=
			\sum_{j=2}^d
				(2-v_j)^{-d\de}.
	$$
$f$ is convex (since its summands are convex), so the maximum value it takes on $T_1^j(\De)\supset R_k\cap\De_1 $ is obtained at one of its $d$ vertices. By symmetry, this maximum is either $f(e_1)$ or $f\big(\frac k{k+1} e_1 + \frac 1{k+1} e_j\big)$, and we find it is the latter:
	$$
			f(e_1) 
		= 
			2^{-d\de}(d-1) 
		\leq 
        	\left(
        		\frac
        			{k+1}
        			{2k+1}
        	\right)^{d\de} 
		+ 
        	2^{-d\de}
        	(d-2)
        =
   			f\left(
   				\frac k{k+1} e_1 + \frac 1{k+1} e_j
   			\right),
	$$
since $\frac k{k+1} e_1 + \frac 1{k+1} e_2$ lies in $R_k\cap\De_1$, $a_k$ takes the claimed value.

%		\begin{figure}
%					\centering
%					\begin{tikzpicture}[scale=8]
%		
%					 	 \foreach \n/\r in {3/1.732}
%					 	 {
%							\path[fill=yellow]
%								(1/\n,0) -- 
%								(0.5/\n,0.5*\r/\n) -- 
%								({0.75/(\n-1)}, {0.25 * \r/(\n-1)}) -- 
%								(1/\n,0);
%								
%							\draw
%								(1/\n,0) --
%								(0.5/\n,0.5*\r/\n) --
%								({0.5/(\n-1)},{0.5*\r/(\n-1)}) --
%								({1/(\n-1)},0) --
%								(1/\n,0) node[midway,anchor=north]{$\qquad R_k \cap \De_1$};
%					 	 	
%					 	 	\draw[fill=black] (1/\n,0) 
%					 	 		node [anchor = east] 
%					 	 			{$\left(
%					 	 				\frac {k+1}{k+2} , 0, \frac 1{k+2} 
%					 	 			\right)$}
%					 	 		circle (.005) ;
%					 	 	
%					 	 	\draw[fill=black] ({0.5/(\n-1)},{0.5*\r/(\n-1)})%({0.75/(\n-1)}, {0.25 * \r/(\n-1)})
%					 	 		node [anchor = west] 
%					 	 			{$\left(
%					 	 				 \frac {k}{k+1}, \frac 1{k+1}, 0
%					 	 			\right)$}
%					 	 		circle (.005) ;	 	 	
%					 	 }
%					 \end{tikzpicture}
%					 \caption[Two maximisers on $R_k\cap \De_1$]{Two maximisers giving the values of $a_k$ and $b_k$ in Lemma 4.4, when $d=3$.}
%					 \label{figr-maximisers on R_k}
%		\end{figure}
%		\label{lemr-values of a_k, b_k}

\medskip

  {\footnotesize

  \noindent
  \textsc{M. Pollicott, Mathematics Department, University of Warwick, Coventry, CV4 7AL, UK.}

  \par

  \addvspace{\medskipamount}

  \noindent
  \textsc{B. Sewell, Alfr\'ed R\'enyi Institute, 13--15 Re\'altonoda utca, Budapest 1053, Hungary.}

}

\end{document}